\def\timestamp{%
Time-stamp: <connected-F-space.tex: Thursday 11-03-2004 at 12:00:00 (cet)>}
\def\stripname Time-stamp: <#1 #2>{#2}
\edef\filedate{\expandafter\stripname\timestamp}
\documentclass[a4paper]{amsart}
\usepackage[numeric]{amsrefs}[2002/02/28]

\newcommand{\bS}{\mathbf{S}}
\newcommand{\cl}{\operatorname{cl}}
\newcommand{\Fr}{\operatorname{Fr}}
\newcommand{\Int}{\operatorname{Int}}
\DeclareMathSymbol\restr{2}{AMSa}{"16}
\DeclareMathSymbol\R{0}{AMSb}{`R}
\newcommand{\cont}{\mathfrak{c}}
\newcommand{\1}{\mathfrak{1}}

\theoremstyle{plain}
\newtheorem{lemma}{Lemma}[section]
\newtheorem{proposition}[lemma]{Proposition}

\theoremstyle{definition}
\newtheorem{question}{Question}

\begin{document}
\title{A connected $F$-space}
\author{Klaas Pieter Hart}
\address{Faculty of Electrical Engineering, Mathematics, and
         Computer Science\\
         TU Delft\\
         Postbus 5031\\
         2600~GA {} Delft\\
         the Netherlands}
\email{K.P.Hart@EWI.TUDelft.NL}
\urladdr{http://aw.twi.tudelft.nl/\~{}hart}
\date{\filedate}

\begin{abstract}
We present an example of a compact connected $F$-space with a continuous 
real-valued function $f$ for which the set 
$\Omega_f =\bigcup\{\Int f^\gets(x):x\in\R\}$
is not dense. 
This indirectly answers a question from Abramovich and Kitover in 
the negative.
\end{abstract}

\keywords{continuum, \textit{F}-space, $d$-independence, $d$-basis}
\subjclass[2000]{Primary: 54G20. Secondary: 46A40 54D30  54F15 54G05}

\maketitle

\section*{Introduction}

The purpose of this note is to give a positive answer to 
Problem~4 from~\cite{Ab-Ki2003}.
The problem asks whether there are a compact and connected $F$-space~$K$
and a continuous real-valued function~$f$ on~$K$ such that
the set~$\Omega_f$ is not dense in~$K$,
where $\Omega_f=\bigcup\{\Int f^\gets(x):x\in\R\}$.
If $K$~is such a space then the vector lattice $C(K)$ has a maximal 
$d$-independent system that is not a $d$-base, which answers Problem~1
from the same paper in the negative.

As defined in~\cite{Ab-Ki2003} a \emph{$d$-independent system} in a vector 
lattice~$X$ is a subset~$D$ with the property that for every band~$B$ in~$X$, 
for every finite subset~$F$ of~$D$ and every choice $\{c_d:d\in F\}$ of 
non-zero scalars the condition $\sum_{d\in F}c_dd\perp B$ implies 
$d\perp B$ for all $d\in F$.
A $d$-independent system~$D$ is a \emph{$d$-basis} if for every $x\in X$
one can find a full system $\mathcal{B}$ of pairwise disjoint bands and
a subset~$\{y_B:B\in\mathcal{B}\}$ of~$X$ such that for each~$B$ the
element~$y_B$ is a linear combination of members of~$D$ and $x-y_B\perp B$.

In topological terms a $d$-independent system in~$C(K)$ is a subset~$D$
such that for every nonempty open subset~$O$ the family of nonzero
members of $\{d\restr O:d\in D\}$ is linearly independent. 
The $d$-independent set~$D$ is a $d$-basis if for each $g\in C(K)$ 
there is a pairwise disjoint family~$\mathcal{O}$ of open sets with a dense
union and such that for every $g\in C(K)$ and every $O\in\mathcal{O}$ the 
restriction~$g\restr O$ is a linear combination of finitely members of 
$\{d\restr O:d\in D\}$.

As observed in~\cite{Ab-Ki2003} for our example $K$ the set $\{\1\}$, 
consisting of just the constant
function with value~$1$, is maximally $d$-indepent in~$C(K)$.
Indeed, if $g$~is not constant then its image~$g[K]$ is a nontrivial
interval; we let $t$ be its mid-point.
Because $K$ is an $F$-space the closed sets 
$\cl g^\gets\bigl[(-\infty,t)\bigr]$ and 
$\cl g^\gets\bigl[(t, \infty)\bigr]$ are disjoint
and because $K$~is connected they do not cover~$K$.
The nonempty open set $\Int g^\gets(t)$ now witnesses that $\{\1,g\}$
is not $d$-independent.
The continuous function~$f$, on the other hand, witnesses that $\{\1\}$
is not a $d$-basis, for clearly any `$d$-linear combination'~$g$ 
of~$\{\1\}$ must have its set~$\Omega_g$ dense in~$K$.

\section{The example}

Let $S$ be the unit square, i.e., $S=[0,1]^2$.
We consider the product $\bS=\omega\times S$, its \v{C}ech-Stone
compactification $\beta\bS$ and the extension $\beta\pi$ of the 
map $\pi:\bS\to\omega$, defined by $\pi(n,x)=n$.

For each free ultrafilter $u\in\beta\omega\setminus\omega$ the fiber
$S_u=\beta\pi^\gets(u)$ is a continuum --- see, e.g., \cite{Ha}.
As it is a closed subset of the \v{C}ech-Stone remainder $\bS^*$ it
is also a compact $F$-space.

The function $f:\bS\to[0,1]$, defined by $f(n,x,y)=x$ is clearly
continuous; we write $f_u$ for the restriction of $\beta f$ to $S_u$.
We shall find a continuum~$K$ in~$S_u$ such that $g=f_u\restr K$ is as
required, i.e., $\Omega_g$~is not dense in~$K$.

We need to describe the boundaries of the fibers of~$f$.
We define $L_t=f_u^\gets(t)\cap \cl f_u^\gets\bigl[[0,t)\bigr]$
and $R_t=f_u^\gets(t)\cap\cl f_u^\gets\bigl[(t,1]\bigr]$; note that
$L_0=R_1=\emptyset$. 

\begin{lemma}\label{lemma.two.components}
For each $t\in(0,1)$ the sets $L_t$ and $R_t$ are exactly the 
components of the boundary~$\Fr f_u^\gets(t)$ of $f_u^\gets(t)$.
\end{lemma}
\begin{proof}
Because $S_u$ is an $F$-space the closed sets $L_t$ and $R_t$ are disjoint;
they cover $\Fr f_u^\gets(t)$ and, because $S_u$ is connected, both are
nonempty.
This shows that $\Fr f_u^\gets(t)$ has at least two components.

To finish we show that $L_t$ and $R_t$ are connected.
For this we first observe that the `rectangle'
$P_{s,r}=S_u\cap\cl\bigl(\omega\times[s,r]\times[0,1]\bigr)$
is connected whenever $s<r$.
This in turn implies that $L_{s,t}=\cl\bigcup_{s<r<t}P_{s,r}$ is connected
whenever $s<t$.
It is readily verified that $L_t=\bigcap_{s<t}L_{s,t}$, hence $L_t$ is 
connected as the intersection of a chain of continua.
By symmetry $R_t$ is also connected.
\end{proof}

This argument also shows that $R_0=\Fr f_u^\gets(0)$ and 
$L_1=\Fr f_u^\gets(1)$ are connected.

We need some more notation.
We denote by $B_u$ the intersection of~$S_u$
with the closure, in $\beta\bS$, of $\omega\times[0,1]\times\{0\}$ 
--- the bottom line of~$S_u$ ---
and likewise  the top line $T_u$ is 
$S_u\cap\cl_{\beta\bS}(\omega\times[0,1]\times\{1\})$.
The continuum $K$ will be defined as the union of the bottom line
of $S_u$ and a family of vertical continua, each of which meet both
the bottom and top lines.

To define this family we define sequences 
$\langle X_\alpha\rangle_\alpha$ and 
$\langle f_\alpha\rangle_\alpha$ 
of closed sets and functions respectively, by recursion.
To begin let $X_0=S_u$.
Given $X_\alpha$ put $f_\alpha=f_u\restr X_\alpha$ and define
$X_{\alpha+1}=X_\alpha\setminus\bigcup_t\Int_\alpha f_\alpha^\gets(t)$,
where $\Int_\alpha$ is the interior operator in $X_\alpha$.
If $\alpha$ is a limit we just let $X_\alpha=\bigcap_{\beta<\alpha}X_\beta$.

\begin{lemma}
For every $\alpha$ and every $t$ the intersections $X_\alpha\cap L_t$ 
and $X_\alpha\cap R_t$ are nonempty
\end{lemma}

\begin{proof}
The proof is by induction on~$\alpha$.

The statement is clearly true for $\alpha=0$ and the case $\alpha=1$ is covered
by Lemma~\ref{lemma.two.components}, whose proof also establishes the
successor step in the induction.
Indeed, to show that $X_{\alpha+1}\cap L_t\neq\emptyset$ we note that, by the 
inductive assumption we know that $P_{s,r}\cap X_\alpha$ meets $L_q$ and 
$R_q$, whenever $s<q<r$.
Therefore $L_{s,t}\cap X_\alpha\neq\emptyset$ for all $s<t$; using compactness
we find that $L_t\cap X_{\alpha+1}=\bigcap_{s<t}(L_{s,t}\cap X_\alpha)$
is nonempty.

The case of limit $\alpha$ follows using compactness as well.
\end{proof}

\begin{lemma}
Every component of $X_\alpha$ meets both $B_u$ and $T_u$.  
\end{lemma}

\begin{proof}
This is clear when $\alpha=0$ and as in the previous lemma
we draw inspiration from the proof of Lemma~\ref{lemma.two.components} for
the argument in the successor step.
Observe first that a component of $X_{\alpha+1}$ is necessarily a subset
of some $L_t$ or $R_t$: these sets are the components of~$X_1$.

Let $C$ be a component of $L_t$ and let $O$ be an arbitrary clopen 
neighbourhood of~$C$ in~$L_t\cap X_{\alpha+1}$; 
choose open sets $U$ and $V$ in~$S_u$ with disjoint closures such 
that $O\subseteq U$ and 
$(L_t\cap X_{\alpha+1})\setminus O\subseteq V$.
There is an~$s$ such that $L_{s,t}\cap X_\alpha\subseteq U\cup V$.
Choose $r\in (s,t)$ such that some component, $D$, 
of~$X_\alpha\cap(L_r\cup R_r)$ meets~$U$; then $D\subseteq U$ and it follows
that $U$ intersects both $B_u$ and~$T_u$.
Because $O$ and $U$ were arbitrary it follows that $C$~must meet~$B_u$ 
and~$T_u$ as well.

In case $\alpha$ is a limit and $C$ a component we have 
$C=\bigcap_{\beta<\alpha}C_\beta$, where $C_\beta$ is the component 
of~$X_\beta$ that contains~$C$; the $C_\beta$'s form a chain and all of them
intersect $B_u$ and $T_u$ and hence by compactness so does $C$.
\end{proof}

There will be a minimal ordinal $\delta$ such that $X_\delta=X_{\delta+1}$
(some information on $\delta$ will be given in the next section).
This means that $\Int_\delta f_\delta^\gets(t)=\emptyset$ for all~$t$.

Our continuum~$K$ is the union of~$B_u$ and $X_\delta$.
Because all components of $X_\delta$ meet $B_u$ we know that $K$~is indeed
connected.
Because each component meets $T_u$ we know that $K$ reaches all the way
up to~$T_u$; by the choice of $\delta$ we get that 
$\Int_Kg^\gets(t)\subseteq B_u$ for all~$t$.
Thus $\Omega_g\subseteq B_u$ and the latter set is certainly not dense in $K$.

\section{A Remark and a question}

The first (and erroneous) version of $K$ was simply 
$B_u\cup\bigcup_{0<t\le1}R_t\cup\bigcup_{0\le t<1}L_t$.
After I realized that the restriction of~$f$ to this subspace did not
provide an example it became clear that the procedure of removing interiors
of fibers had to be iterated, which lead to the sequence
$\langle X_\alpha\rangle_\alpha$.
We can provide some information on the ordinal~$\delta$ at which the sequence
becomes constant.

\begin{proposition}
$\delta<\cont^+$  
\end{proposition}

\begin{proof}
Let $\mathcal{B}$ be a base for $S_u$ of cardinality $\cont$.
For every $\alpha<\delta$ there is a $B_\alpha\in\mathcal{B}$ such that
$\emptyset\neq B_\alpha\cap X_\alpha\subseteq X_\alpha\setminus X_{\alpha+1}$.
Clearly $\alpha\mapsto B_\alpha$ is one-to-one, which establishes that 
$|\delta|\le\cont$.
\end{proof}

The $F$-space property implies that $\delta$ cannot be a successor ordinal,
nor an ordinal of countable cofinality.

\begin{lemma}\label{lemma.homogeneous}
If $\alpha<\delta$ then $X_\alpha\setminus X_{\alpha+1}$ meets every
$L_t$ and every $R_t$.  
\end{lemma}

\begin{proof}
This is basically a consequence of the homogeneity of the unit interval.
If $h:[0,1]\to[0,1]$ is a homeomorphism then it induces an 
autohomeomorphism~$h_u$ of~$S_u$ via the map $(n,x,y)\mapsto(n,h(x),y)$
from~$\bS$ to itself.
The map~$h_u$ simply permutes the fibers~$f^\gets(t)$ and it is relatively
straightforward to show by induction that $h_u[X_\alpha]=X_\alpha$
for all~$\alpha$. 
There are enough maps~$h$ to ensure that once 
$X_\alpha\setminus X_{\alpha+1}$ meets one~$L_t$ (or one~$R_t$) it meets 
all~$L_s$ and all~$R_s$.
\end{proof}

\begin{proposition}
$\delta$ is not a successor ordinal.  
\end{proposition}

\begin{proof}
Let $\alpha<\delta$, we show that $\alpha+1<\delta$.
Fix $t\in(0,1)$ and let $\langle t_n\rangle_n$ be a sequence in~$[0,1]$
that converges to~$t$ from above.
By Lemma~\ref{lemma.homogeneous} we can pick 
$x_n\in L_{t_n}\cap X_\alpha\setminus X_{\alpha+1}$ for each~$n$.

Clearly every point in the closure of $\{x_n\}_n$ belongs to $X_{\alpha+1}$;
we show that none belong to~$X_{\alpha+2}$.
To see this observe that the $F_\sigma$-sets $F=\{x_n\}_n$ and 
$G=f^\gets\bigl[(t,1]\bigr]$ are \emph{separated} in~$S_u$,
i.e., $\cl F\cap G=\emptyset=F\cap\cl G$.
Using normality in the form of Urysohn's Lemma one can find a continuous
function $h:S_u\to[-1,1]$ such that $h[F]\subseteq[-1,0)$ and
$h[G]\subseteq(0,1]$.
But now the $F$-space property applies to show that 
$\cl F\cap\cl G=\emptyset$.
\end{proof}

In a similar way we can prove the following.

\begin{proposition}
The ordinal $\delta$ has uncountable cofinality.  
\end{proposition}

\begin{proof}
We choose an increasing sequence $\langle\alpha_n\rangle_n$
of ordinals below~$\delta$; we show that $\lim_n\alpha_n<\delta$.

As in the previous proof we fix $t\in(0,1)$ and a 
sequence $\langle t_n\rangle_n$ converging to~$t$ from above.
As before we choose $x_n\in L_{t_n}\cap X_{\alpha_n}\setminus X_{\alpha_n+1}$
for all~$n$.

As in the previous proof the $F$-space property now ensures that every point 
in the closure of~$\{x_n\}_n$
belongs to~$X_\alpha\setminus X_{\alpha+1}$.
\end{proof}

We deduce that $\delta$ must be at least $\omega_1$ but the following
question remains.

\begin{question}
What is the exact value of $\delta$?  
\end{question}


\begin{bibdiv}
\begin{biblist}

\bib{Ab-Ki2003}{article}{
    author={Abramovich, Y. A.},
    author={Kitover, A. K.},
     title={$d$-Independence and $d$-bases},
   journal={Positivity},
    volume={7},
      date={2003},
     pages={95\ndash 97},
      issn={1385-1292},
}



\bib{Ha}{incollection}{
    author={Hart, Klaas~Pieter},
     title={The \v{C}ech-Stone compactification of the Real Line},
      date={1992},
 booktitle={Recent progress in general topology},
    editor={Hu\v{s}ek, Miroslav},
    editor={van Mill, Jan},
 publisher={North-Holland},
   address={Amsterdam},
     pages={317\ndash 352},
    review={MR 95g:54004},
}

\end{biblist}
\end{bibdiv}
\end{document}